\input amstex
\documentstyle{amsppt}
\input bull-ppt
\define\GZ{\operatorname{ GL}_{2}(\Bbb Z)}
\define\GR{\operatorname{ GL}_{2}(\Bbb R)}
\define\DFAF{|D_{F}|^{1/n(n-1)}A_{F}}

\let\sl\it
\topmatter
\cvol{31}
\cvolyear{1994}
\cmonth{October}
\cyear{1994}
\cvolno{2}
\cpgs{204-207}
\title An isoperimetric inequality related to Thue's 
equation\endtitle
\shorttitle{ 
An Isoperimetric Inequality}
\author Michael A.\ Bean\endauthor
\address Mathematical Sciences Research Institute, 1000 
Centennial Drive,
Berkeley, California 94720\endaddress 
\cu 
Department of Mathematics, University of Toronto, Toronto,
Ontario  M5S~1A1, Canada\endcu 
\date August 27, 1993\enddate
\ml 
bean\@math.toronto.edu\endml 
%
\subjclass Primary 11D75, 11J25, 51M25; Secondary 26D20, 
33B15\endsubjclass 
\abstract This paper announces the discovery of an 
isoperimetric
inequality for the area of plane regions defined by binary 
forms.
This result has been applied subsequently in the 
enumeration of solutions to the
Thue inequality and, given its fundamental nature,
 may find application in other areas as well. 
\endabstract 
\thanks Supported in part by a Postdoctoral Fellowship 
from the 
Natural Sciences and Engineering Research Council of 
Canada\endthanks 
\endtopmatter

\document

\heading 1. Introduction \endheading
Let $F(X,Y) = a_{0}X^{n} + a_{1}X^{n-1}Y + \cdots + 
a_{n}Y^{n}$
be a binary form with rational integer coefficients, and 
let $A_{F}$
denote the area of the region $|F(x,y)| \leq 1$ in the 
real affine plane.
In a seminal paper of 1909, Thue \cite{11} showed that if 
$F$ is an
irreducible form with degree at least three and if $h$ is 
a non-zero
integer, then the equation $F(x,y) = h$ has only a finite 
number of 
solutions in integers $x$~and~$y$.
In~1934, Mahler \cite{7} gave an estimate for the number, 
$N_{F}(h)$, of
solutions of the Thue inequality $|F(x,y)| \leq h$ in terms
of the area, $A_{F}h^{2/n}$, of this region; in 
particular, he showed,
under the same conditions as Thue, that
$$ | N_{F}(h) - A_{F}h^{2/n} | \leq c h^{1/(n-1)} $$
where $c$ is a number which depends only on $F$. 
Mahler left $c$ and $A_{F}$ unspecified, but did show that 
$A_{F}$
is finite (under the above conditions).

Surprisingly little attention has been given to the 
estimation of $A_{F}$,
even though it is clear that sharp estimates for $A_{F}$ 
would have
significant implications for Thue inequalities. Mueller 
and Schmidt
\cite{8,9}
showed, under certain technical conditions, that $A_{F}$ 
is bounded when
the number of coefficients of $F$ is small, but their 
approach gave
little insight into the nature of $A_{F}$ in general. The 
purpose of 
this paper is to announce the discovery of an isoperimetric
inequality for $A_{F}$. Detailed proofs will be 
forthcoming in \cite{4}.

\proclaim{Theorem 1} Let $F$ be a binary form with integer 
coefficients
having non-zero discriminant and degree at least three. Then
$$ A_{F} \leq 3 B\left(\frac{1}{3},\frac{1}{3}\right) $$
where $B(\frac{1}{3},\frac{1}{3})$ denotes the Beta 
function with 
arguments of $\frac13$. This bound is attained for forms 
with integer
coefficients which are equivalent under $\GZ$ to $XY(X-Y)$. 
\endproclaim
The approximate numerical value of $3 
B(\frac{1}{3},\frac{1}{3})$ is
$15.8997$.
A similar result is true for forms with complex 
coefficients provided that 
one includes a factor involving the discriminant of $F$
(see \S 2). However, the conditions on the discriminant
and the degree cannot be relaxed in general since, for 
example,
the forms $X^{n}$ and $X^{2}-Y^{2}$ give rise to infinite 
area.
The condition that $F$ have non-zero discriminant is 
equivalent to 
the condition that neither $F(x,1)$ nor $F(1,y)$ have 
multiple roots
and is satisfied, in particular, when $F$ is an 
irreducible form.
The notion of equivalence under $\GZ$ will be explained in 
\S 2.

One of the striking aspects of Theorem~1 is that the bound 
given
for $A_{F}$ is relatively small. In light of Mahler's 
result,
one might expect $N_{F}(1)$ to be bounded and small as well.
However, an examination of the forms
$$ P_{k}(X,Y) = X^{2k} + (X-Y)^{2}(2X-Y)^{2} \cdots 
(kX-Y)^{2} $$
reveals that this is not the case. Indeed, the region 
$|P_{k}(x,y)| \leq 1$
contains at least $2k$ integer lattice points, namely,
$$ \pm (1,1), \pm (1,2), \ldots, \pm (1,k), $$
but nevertheless is bounded (since $P_{k}$ is positive 
definite)
and has area less than~$16$. (See Figure~1 for the graph 
of $P_{7}$.)

These examples suggest that a finer analysis of the error 
term in Mahler's
result could be quite delicate. However, they raise the 
interesting
question of classifying forms according to the size of 
this error.


\fg{21pc}
\caption{Figure 1. $|P_7(x,y)|=1$}
\endfg

\heading 2. The main ideas \endheading
There are three main ideas involved in the proof of 
Theorem~1.
The first and most important is to consider the estimation
of $A_{F}$ over the larger class of forms having complex 
coefficients.

For any form $F$, let $D_{F}$ denote the discriminant of 
$F$.
If $F$ has a factorization 
$\prod_{i=1}^{n}(\alpha_{i}X-\beta_{i}Y)$
with $\alpha_{i}$, $\beta_{i} \in \Bbb C$ (every binary form
with complex coefficients has such a factorization), then
$D_{F} = \prod_{i<j}(\alpha_{i}\beta_{j}-\alpha_{j}%
\beta_{i})^{2}$.
Let $\GR$ denote the group of $2 \times 2$ real invertible 
matrices.
For each $T=\txtpmatrix a & b \\  c & d \endtxtpmatrix \in 
\GR$, let
$F_{T}(X,Y)=F(aX+bY,cX+dY)$. Two forms are said to be 
equivalent
under $\GR$ if $G=F_{T}$ for some $T \in \GR$.
Equivalence under $\GZ$ is defined similarly.

In this more general setting, we can prove the following 
inequality.
\proclaim{Theorem 2} Let $F$ be a binary form with complex 
coefficients
having degree $n \geq 3$ and discriminant $D_{F} \neq 0$. 
Then
$$ \DFAF \leq 3 B \left( \frac{1}{3},\frac{1}{3}\right). $$
This bound is attained precisely when $F$ is a cubic form 
which,
up to multiplication by a complex number, is equivalent 
under $\GR$ 
to the form $XY(X-Y)$.
\endproclaim
Since the discriminant of a form with integer coefficients 
is an
integer, Theorem~1 is an immediate consequence.
The inequality in Theorem~2 is an {\sl isoperimetric\/} 
inequality
since it is sharp and connects the geometric quantities
$A_{F}$ and $|D_{F}|$ which are related to the shape and 
size of the region $|F(x,y)| \leq 1$.
The quantity
$\DFAF$ is natural to consider since it is invariant under 
$\GR$
(i.e., $|D_{F_{T}}|^{1/n(n-1)}A_{F_{T}} = \DFAF$ for all 
$T \in \GR$)
in the same way that $A_{F}$ is invariant under $\GZ$.

The real advantage of considering $\DFAF$ instead of $A_{F}$
is that one can apply the inequality between arithmetic 
and geometric
means in a systematic fashion to reduce the estimation of 
$\DFAF$
to cubic forms. This is the second main idea. To be 
specific,
let
$$ M_{n} = \max \DFAF $$
where the maximum is taken over all forms of degree $n$ 
with 
$D_{F} \neq 0$. Then one can show that $\{M_{n}\}$ is a 
decreasing 
sequence.

The third main idea is to consider the quantity $\DFAF$ as 
a function
of $n$ complex variables (the roots of $F(x,1)$).
By appealing to an appropriate maximum principle, one can 
then deduce
that $\DFAF$ is maximized by a form $F$ for which $F(x,1)$ 
has $n$
distinct real roots. Heuristically, one expects $\DFAF$ to 
be
maximized by a form whose corresponding graph $|F(x,y)| = 
1$ has
the maximal number of asymptotes in the real affine plane
(i.e.\ by a form for which $F(x,1)$ has $n$~distinct real 
roots).
In the case of cubics, every such form is equivalent under 
$\GR$
to $XY(X-Y)$.

\heading 3. Related questions \endheading
The outline given above raises some interesting questions 
regarding
the sequence~$\{M_{n}\}$. In particular, one might like to 
know
all the values of this sequence and its limiting behavior.
These and other questions will be addressed in \cite{3}.
However, based on that work, we believe that the following 
is true.
\dfn 
{Conjecture}
The maximal value $M_{n}$ of the quantity $\DFAF$
over the class of forms of degree $n$ with complex 
coefficients
and non-zero discriminant is attained precisely when $F$ 
is a form
which, up to multiplication by a complex number, is 
equivalent
under $\GR$
to the form
$$ F_{n}^{*}(X,Y) = \prod_{k=1}^{n} \left(X \sin 
\left(\frac{k\pi}{n}\right) -
                         Y \cos \left( \frac{k\pi}{n} 
\right) \right) .$$
Moreover, 
the sequence of maximal values $\{M_{n}\}$ 
decreases monotonically
to the value $2\pi$.
\enddfn 

\heading Acknowledgment \endheading
The author is grateful to E.~Bombieri, J.~Friedlander, 
C.~Stewart, and the referee 
for their suggestions and comments.

\enddocument